\documentclass[12pt]{article}
\usepackage{amsfonts}
\usepackage[dvips]{graphics}
\newtheorem{theorem}{Theorem}

\newtheorem{definition}{Definition}
\newtheorem{remark}{Remark}
\newtheorem{example}{Example}
\newtheorem{corollary}{Corollary}
\makeindex
\begin{document}
\title{A  Liouville comparison principle  for solutions of semilinear parabolic
second-order partial differential inequalities.}
\author{Vasilii V. Kurta.}
\maketitle
\thispagestyle{empty}
\begin{abstract}
\noindent \noindent We obtain a new Liouville comparison principle
for entire weak solutions $(u,v)$ of semilinear parabolic
second-order partial differential inequalities of the form
$$
\qquad \qquad \qquad u_t -{\mathcal L}u- |u|^{q-1}u\geq v_t -{\mathcal L}v- |v|^{q-1}v\qquad \qquad \qquad   (*)
$$
in the half-space  ${\mathbb S} = {\mathbb R}^1_+ \times \mathbb
R^n$. Here $n\geq 1$, $q>0$ and
$$
{\mathcal L}=\sum\limits_{i,j=1}^n\frac{\partial
}{{\partial}x_i}\left [ a_{ij}(t, x) \frac{\partial }{{\partial
}x_j}\right],$$ where   $a_{ij}(t,x)$, $i,j=1,\dots ,n$, are
functions defined, measurable and locally bounded in $\mathbb S$,
and such that $a_{ij}(t,x)=a_{ji}(t,x)$ and
$$ \sum_{i,j=1}^n a_{ij}(t,x)\xi_i\xi_j\geq 0
$$
for almost all $(t,x)\in \mathbb S$ and all $\xi \in \mathbb R^n$.
The critical exponents in the Liouville comparison principle
obtained, which responsible for the non-existence of non-trivial
(i.e., such that $u\not \equiv v$) entire weak solutions to ($*$) in
$\mathbb S$, depend on the behaviour of the coefficients of the
operator $\mathcal L$  at infinity. As direct corollaries we obtain
a new Fujita comparison principle for entire  weak solutions $(u,v)$
of  the Cauchy problem for the inequality ($*$), as well as new
Liouville-type and Fujita-type theorems for non-negative entire weak
solutions $u$ of the inequality ($*$) in  the case when $v\equiv 0$.
All the results obtained are new and sharp. \end{abstract}

\section{Introduction and preliminaries}

This work is devoted  to a new Liouville comparison principle  of
elliptic type for entire weak solutions to parabolic inequalities of
the form
\begin{eqnarray}
 u_t -{\mathcal L}u - |u|^{q-1}u\geq v_t -{\mathcal L}v - |v|^{q-1}v
\end{eqnarray}
in the half-space  ${\mathbb S} = (0, +\infty) \times \mathbb R^n$,
where $n\geq 1$ is a natural number,  $q>0$ is a real number and
${\mathcal L}$ is a linear  second-order  partial differential
operator  in divergence form defined   by the relation
\begin{eqnarray}
{\mathcal L}=\sum\limits_{i,j=1}^n\frac{\partial }{{\partial
}x_i}\left [ a_{ij}(t, x) \frac{\partial }{{\partial }x_j}\right].
\end{eqnarray}
Here and in what follows, we assume that  the coefficients
$a_{ij}(t, x)$, $i,j=1, \dots, n$, of the operator ${\mathcal L}$
are functions defined, measurable and locally bounded  in ${\mathbb
S}$, and such that $a_{ij}(t, x)=a_{ji}(t, x)$, $i,j=1,\dots ,n$,
for almost all $(t, x)\in {\mathbb S}$. Also,  we assume that  the
corresponding quadratic form satisfies the conditions
\begin{equation} 0\leq \sum\limits_{i,j=1}^n
a_{ij}(t,x)\xi_i\xi_j\leq A(t, x)|\xi|^2 \end{equation} for all
$\xi=(\xi_1, \dots, \xi_n)\in {\mathbb R}^n$ and  almost all $(t,
x)\in {\mathbb S}$, with $A(t, x)$ some function defined,
measurable, non-negative  and locally bounded in ${\mathbb S}$.

It is important to note that if $u=u(t,x)$ and $v=v(t,x)$  satisfy
inequalities \begin{eqnarray} u_t \geq {\mathcal L}u+|u|^{q-1}u
\end{eqnarray} and  \begin{eqnarray} v_t \leq {\mathcal
L}v+|v|^{q-1}v,\end{eqnarray} then the pair $(u,v)$ satisfies the
inequality (1). Thus, all the results obtained in this paper for
solutions of (1) are valid for the corresponding solutions of the
system (4)--(5).

Under entire solutions of  inequalities (1), (4) and  (5)   we
understand solutions defined in the whole half-space $\mathbb S$,
and under Liouville results of elliptic type for solutions of
evolution inequalities (1), (4) and  (5)  in the half-space $\Bbb S$
we understand Liouville-type results which, in their formulations,
have no restrictions on the behaviour of solutions to these
inequalities on the hyper-plane $t=0$.  Also, we would like to
underline that we impose neither growth conditions on the behaviour
of solutions to inequalities (1), (4) and  (5) or on that of any of
their partial derivatives at infinity.

In the case when the coefficients of the operator $\mathcal L$ are
globally bounded in $\mathbb S$,   a Liouville comparison principle
of elliptic type for entire weak solutions $(u,v)$ of the inequality
(1), as well as  Liouville-type and Fujita-type theorems for
non-negative entire weak solutions $u$ of the inequality (4), were
obtained in [7]. In those results, a critical exponent which is
responsible for the non-existence of non-trivial (i.e., such that
$u\not\equiv v$) entire weak solutions $(u,v)$  to the inequality
(1), as well as  non-trivial (i.e., such that $u\not\equiv 0$)
non-negative entire weak solutions $u$ to the inequality (4),
coincides with the well-known Fujita critical blow-up exponent for
non-trivial non-negative entire classical solutions to the Cauchy
problem for the equation
\begin{eqnarray}
u_t -\Delta u =    |u|^{q-1}u,
\end{eqnarray}
which was established  in [2], [4] and [8]. However, it is
intuitively clear that the character of the behaviour of the
coefficients $a_{ij}(t,x)$ of the operator $\mathcal L$ as $|x|\to +
\infty$  must manifest itself in Liouville-type and  Fujita-type
results. In particular, a potential critical exponent in a Liouville
comparison principle for  entire weak solutions of  (1), which is
responsible for the non-existence of non-trivial entire weak
solutions to the inequality (1) must depend on the behaviour of the
coefficients of the operator $\mathcal L$ as $|x|\to + \infty$.

In order to trace this dependence we consider the value
\begin{eqnarray}{\mathcal A}(R)={\mathrm {ess}\sup}_{(t,x)\in
(0,+\infty )\times \{R/2 <|x|<R \}}A(t,x)\end{eqnarray} for any
$R>0$ and assume that the coefficients of the operator $\mathcal L$
satisfy  the condition
\begin{eqnarray}{\mathcal A}(R)\leq cR^{2-\alpha},\end{eqnarray}
with some real constant $\alpha$ and some real positive constant
$c$, for all $R>1$. It is clear that if $\alpha <2$, then the
coefficients of the operator $\mathcal L$ may be unbounded in
$\mathbb S$, if $\alpha =2$, the coefficients of the operator
$\mathcal L$ are globally bounded in $\mathbb S$, and if $\alpha
>2$, they must vanish as $|x|\to +\infty$. Our main concern in this
paper  is the cases when $\alpha\neq 2$.

We also introduce a special function space, which is directly
associated to the linear partial differential operator
$\displaystyle {\mathcal P}= \frac{\partial}{\partial t}- {\mathcal
L}, $ and assume that entire weak solutions of  inequalities (1),
(4) and (5)  belong to this function space only locally  in $\mathbb
S$.

\section{Definitions}

\begin{definition}
Let   $n\geq 1$, $q>0$ and $\hat q =\max\{1,q\}$,   let $\mathcal L$
be a differential operator defined  by (2), and let $\Pi$ be an
arbitrary bounded domain in ${\mathbb S}$. By $W^{\mathcal L,
q}(\Pi)$  we denote the completion of the function space
$C^{\infty}(\Pi)$  with respect to the norm
$$
\| w \|_{W^{\mathcal L,
q}(\Pi)}=
\int\limits_{\Pi}
|w_t|dtdx+
\left[\int\limits_{\Pi}\sum\limits_{i,j=1}^n
a_{ij}(t,x)\frac{\partial w}{\partial x_i}\frac{\partial w}{\partial x_j}dtdx\right]^{1/2} +
\left[\int\limits_{\Pi}|w|^{\hat q}dtdx\right]^{1/{\hat q}},
$$
where $C^{\infty}(\Pi)$ is the space of all functions  defined and
infinitely differentiable  in  $\Pi$.
\end{definition}
\begin{definition}
Let  $n\geq 1$ and  $q>0$, and let $\mathcal L$ be a differential
operator defined  by (2). A function $w=w(t,x)$ belongs to the
function space $W^{\mathcal L, q}_{\mathrm {loc}}({\mathbb S})$ if
$w$ belongs to $W^{\mathcal L, q}(\Pi)$ for any bounded domain $\Pi$
in ${\mathbb S}$.
\end{definition}
\begin{definition}
Let  $n\geq 1$ and  $q>0$, and let $\mathcal L$ be a differential
operator defined  by (2).  A pair $(u,v)$ of functions $u=u(t,x)$
and $v=v(t,x)$ is called  an entire weak  solution to the inequality
(1) in $\mathbb S$, if these functions  are defined and measurable
in ${\mathbb S}$, belong to the function space $W^{\mathcal L,
q}_{\mathrm {loc}}({\mathbb S})$ and satisfy the integral inequality
\begin{eqnarray}
\int\limits_{\mathbb S}\left[u_t\varphi+\sum_{i,j=1}^n
a_{ij}(t,x)\frac{\partial \varphi}{\partial x_i} \frac{\partial
u}{\partial x_j} - |u|^{q-1}u\varphi\right]dtdx \geq \nonumber
\\
\int\limits_{\mathbb S}\left[v_t\varphi+\sum_{i,j=1}^n
a_{ij}(t,x)\frac{\partial \varphi}{\partial x_i} \frac{\partial
v}{\partial x_j}- |v|^{q-1}v\varphi\right]dtdx\end{eqnarray}
for every  function $\varphi \in  C^\infty (\mathbb S)$ with compact
support in $\mathbb S$,  where $C^{\infty}({\mathbb S})$ is the
space of all functions defined and infinitely differentiable  in
${\mathbb S}$.
\end{definition}

\begin{remark} We understand the inequality  (9) in the sense
discussed, e.g., in [10] or [15].
\end{remark}

Analogous definitions of solutions to  inequalities  (4) and  (5),
as special cases of the inequality  (1)  for $v\equiv 0$ or $u\equiv
0$, follow immediately from Definition 3.

\section{Results}
\begin{theorem} Let $n\geq 1$, $\alpha>0$ and  $1<q\leq 1 +\frac \alpha
n$, let $\mathcal L$ be a differential operator defined  by (2), the
coefficients of which  satisfy the condition (8) with the given
$\alpha$ and some $c>0$,   and  let $(u,v)$ be an entire weak
solution of the inequality (1) in $\mathbb S$ such that $u\geq v$.
Then $u = v$ in $\mathbb S$.
\end{theorem}

As we have observed above, since any solutions $u=u(t,x)$,
$v=v(t,x)$ of inequalities (4), (5) is a solution $(u,v)$  of the
inequality (1), then the following statement is a direct corollary
of Theorem 1.
\begin{theorem} Let $n\geq 1$, $\alpha>0$ and  $1<q\leq 1 +\frac \alpha
n$, let $\mathcal L$ be a differential operator defined  by (2), the
coefficients of which   satisfy the condition (8) with the given
$\alpha$ and some $c>0$,   and let $u=u(t,x)$ be an entire weak
solution of the inequality (4) and  $v=v(t,x)$ be an entire weak
solution of the inequality (5) in $\mathbb S$ such that $u\geq v$.
Then $u = v$ in $\mathbb S$.
\end{theorem}

The results in Theorems 1 and 2, which evidently have a comparison
principle character, we term a Liouville-type comparison principle,
since in  particular cases when either $u\equiv 0$ or $v\equiv 0$,
it becomes a Liouville-type theorem for solutions of  inequality (5)
or (4), respectively. We formulate here only  the case when $v\equiv
0$.

\begin{theorem} Let $n\geq 1$, $\alpha>0$ and  $1<q\leq 1 +\frac \alpha
n$, let $\mathcal L$ be a differential operator defined  by (2), the
coefficients of which   satisfy the condition (8) with the given
$\alpha$ and some $c>0$,   and let  $u=u(t,x)$ be a non-negative
entire weak solution of the inequality (4) in $\mathbb S$. Then $u =
0$ in $\mathbb S$.
\end{theorem}

Since in Theorems 1 and 2  we impose no conditions on the behaviour
of entire weak solutions  of   inequalities (1), (4) and   (5)  on
the hyper-plane $t=0$, we can formulate,  as a direct corollary of
the Liouville comparison principle in Theorems 1 and 2, a comparison
principle, which in turn one can term a Fujita comparison principle,
for entire weak solutions of the Cauchy problem with arbitrary
initial data for $u$ and $v$ for  inequalities (1), (4) and  (5) in
$\mathbb S$.

\begin{theorem} Let $n\geq 1$, $\alpha>0$ and $1<q\leq 1 +\frac \alpha
n$, let $\mathcal L$ be a differential operator defined  by (2), the
coefficients of which   satisfy the condition (8) with the given
$\alpha$ and some $c>0$,  and  let $(u,v)$ be an entire weak
solution of the Cauchy problem, with arbitrary initial data for
$u=u(t,x)$ and $v=v(t,x)$,  for the inequality (1) in $\mathbb S$
such that $u\geq v$. Then $u = v$ in $\mathbb S$.
\end{theorem}

Note that  the initial data for $u=u(t,x)$ and $v=v(t,x)$ in Theorem
4 may be different.

\begin{theorem} Let $n\geq 1$, $\alpha>0$ and $1<q\leq 1 +\frac \alpha
n$, let $\mathcal L$ be a differential operator defined  by (2), the
coefficients of which   satisfy the condition (8) with the given
$\alpha$ and some $c>0$,  and let $u=u(t,x)$ be an entire weak
solution of the Cauchy problem,  with arbitrary initial data, for
the inequality (4) and  $v=v(t,x)$ be an entire  weak solution of
the Cauchy problem, with arbitrary initial data,  for the inequality
(5) in $\mathbb S$ such that $u\geq v$. Then $u = v$   in $\mathbb
S$.
\end{theorem}

It is clear that in a particular case when $u\equiv 0$ or $v\equiv
0$, the Fujita comparison principle in Theorems 4 and 5 becomes a
Fujita-type theorem  for entire  weak solutions of the Cauchy
problem for  inequality (5) or (4), respectively. As before, we
formulate here only the case when $v\equiv 0$.

\begin{theorem} Let $n\geq 1$,  $\alpha>0$ and  $1<q\leq 1 +\frac \alpha
n$, let $\mathcal L$ be a differential operator defined  by (2), the
coefficients of  which   satisfy the condition (8) with the given
$\alpha$ and some $c>0$,  and let  $u=u(t,x)$ be a non-negative
entire weak solution of the Cauchy problem, with arbitrary initial
data, for the inequality (4) in $\mathbb S$. Then $u = 0$  in
$\mathbb S$.
\end{theorem}

As we have mentioned above,  if the coefficients of the operator
$\mathcal L$ are globally bounded in $\mathbb S$, then the condition
(8) for these coefficients is fulfilled with $\alpha=2$ and some
constant $c>0$, and, therefore,  the results obtained (here we
restrict ourselves only with Theorem 1) may be formulated in the
following form:
\begin{corollary} Let $n\geq 1$ and  $1<q\leq 1 +\frac 2 n$,
let $\mathcal L$ be a differential operator defined  by (2), the
coefficients of  which  are globally bounded in $\mathbb S$, and let
$(u,v)$ be an entire weak solution of the inequality (1) in $\mathbb
S$ such that $u\geq v$. Then $u = v$ in $\mathbb S$.
\end{corollary}

So, in a  particular case when $\alpha=2$, the critical blow-up
exponent in Theorems 1--6  coincides with the well-known Fujita
critical blow-up exponent, and the well-known Fujita theorem on
blow-up of non-trivial non-negative entire  classical solutions to
the Cauchy problem  with arbitrary initial data  for the equation
(6) proved in [2], [4] and [8] is a direct corollary of Theorem 6
when $\alpha=2$. Also, as we have mentioned above, similar results
to those in Theorems 1--6 when $\alpha=2$ were obtained in [7]. The
difference between the results in Theorems 1--6 when $\alpha=2$ and
those obtained in [7] consists of  the fact that in the present
paper we study  solutions to  inequalities (1), (4) and (5) in the
function space $W^{\mathcal L, q}_{\mathrm {loc}}({\mathbb S})$
which is, generally speaking, wider than that considered  in [7].
Thus, all the results in Theorems 1--6 are new, with new critical
blow-up exponents in the cases when $\alpha\neq 2$.  We demonstrate
their sharpness by the following examples.

\begin{example}
Let $n\geq 1$, $+\infty > \alpha>-\infty $ and $q\leq 1$, and let
$\mathcal L$ be a differential operator defined  by (2), the
coefficients of  which   satisfy the condition (8) with the given
$\alpha$ and some $c>0$. It is evident that the function
$u(t,x)=e^t$ is a positive entire classical solution of the
inequality (4) in $\mathbb S$. Also, it is clear that  the function
$v=-u(t,x)$ is a negative entire classical solution of the
inequality (5) in $\mathbb S$, and, thus, the pair of the  functions
$u=u(t,x)$ and $v=v(t,x)$ is a non-trivial entire classical solution
of the system (4)--(5) and, therefore, of the inequality (1) in
$\mathbb S$ such that $u(t,x)> v(t,x)$.
\end{example}
\begin{example} Let  $n\geq 1$,  $\alpha>0$ and  $q> 1 +\frac \alpha n$.
Consider the operator $\mathcal L$ defined by (2)  with the
coefficients given by the formula
\begin{eqnarray}
a_{ij}(t,x)=(1+|x|^2)^{\frac{2-\alpha}2}\delta_{ij},
\end{eqnarray}
where $\delta_{ij}$ are Kronecker's symbols and $i,j=1,\dots ,n$. It
is easy to see that the condition (8) is fulfilled for these
coefficients with the given $\alpha$ and some $c>0$. Also, for the
given $\alpha$ and $q$, let $\beta=\frac 1{q-1}$, $\frac 1{\alpha
n(q-1)}<\gamma\leq \left( \frac 1{\alpha}\right)^2$, $0<
\kappa\leq\left ( \alpha n\left(\gamma -\frac 1{\alpha
n(q-1)}\right)\right)^{1/(q-1)}$ and
\begin{eqnarray}
u(t,x)=\kappa (t+1)^{-\beta}\exp\left(-\gamma \frac{(1+|x|^2)^{\frac
\alpha 2}}{t+1}\right).\end{eqnarray} Making  necessary
calculations, it is not difficult to verify that the function
$u=u(t,x)$ defined by the formula (11) is a  positive  entire
classical solution of the inequality (4)  in $\mathbb S$, with
$a_{ij}(t,x)$, the coefficients of the operator $\mathcal L$,
defined  by (10). Also, it is clear that the function $v=-u(t,x)$ is
a negative entire classical solution of the inequality (5) in
$\mathbb S$, with $a_{ij}(t,x)$ in (2) defined  by (10), and, thus,
the pair of the functions $u=u(t,x)$ and $v=v(t,x)$ is a non-trivial
entire classical solution of the system (4)--(5) and, therefore, of
the inequality (1) in $\mathbb S$ such that $u(t,x)> v(t,x)$, with
$a_{ij}(t,x)$   in (2) defined by (10).
\end{example}

Note that a  positive entire classical sub-solution of the equation
(6) in a form similar to that given by the formula (11) with
$\alpha=2$ was constructed in [17, p. 283].

\begin{example} Let  $n\geq 1$, $\alpha\leq 0$,  $q>1 +\frac \alpha
n$ and $q>1$,   and let $\hat \alpha$ be any positive number such
that $q>1+\frac {\hat \alpha}n$. Consider the operator $\mathcal L$
defined by (2) with the coefficients given   by the formula
\begin{eqnarray}
a_{ij}(t,x)=(1+|x|^2)^{\frac{2-\hat \alpha}2}\delta_{ij},
\end{eqnarray}
where $\delta_{ij}$ are Kronecker's symbols and $i,j=1,\dots ,n$. As
in Example 2, it is easy to see that $\displaystyle{\mathcal
A}(R)\leq CR^{2-\hat\alpha}$ for all $R>1$, with $C$ some positive
constant which possibly depends on $\hat \alpha$ and $n$, and,
therefore,  the condition (8) is fulfilled for these coefficients
with the given $\alpha$ and some $c>0$. Also, for the given $\hat
\alpha$ and $q$, let $\beta=\frac 1{q-1}$, $\frac 1{\hat \alpha
n(q-1)}<\gamma\leq \left (\frac 1{\hat \alpha}\right)^2$, $0<
\kappa\leq\left ( \hat \alpha n\left(\gamma -\frac 1{\hat \alpha
n(q-1)}\right)\right)^{1/(q-1)}$ and
\begin{eqnarray}
u(t,x)=\kappa (t+1)^{-\beta}\exp\left(-\gamma \frac{(1+|x|^2)^{\frac
{\hat \alpha} 2}}{t+1}\right).\end{eqnarray} Again as in Example 2,
it is not difficult to verify that the function $u=u(t,x)$ defined
by the formula (13) is a  positive  entire classical  solution of
the inequality (4) in $\mathbb S$, with $a_{ij}(t,x)$ in (2) defined
by (12). Also, it is clear that the function $v=-u(t,x)$ is a
negative entire classical solution of the inequality (5) in $\mathbb
S$, with $a_{ij}(t,x)$, the coefficients of the operator $\mathcal
L$, defined by (12), and, thus, the pair of the functions $u=u(t,x)$
and $v=v(t,x)$ is a non-trivial entire classical solution of the
system (4)--(5) and, therefore, of the inequality (1) in $\mathbb S$
such that $u(t,x)> v(t,x)$, with $a_{ij}(t,x)$ in (2) defined by
(12).
\end{example}

\begin{remark} For the case when $n\geq 1$,  $\alpha\leq 0$ and
$1 \geq q> 1 +\frac \alpha
n$, see Example 1.
\end{remark}

Finally, we would like to note  that elliptic analogues of the
results in Theorems 1--6 were obtained in [12] and [13]. To prove
the results obtained in the present work  we further develop an
approach proposed in [11]. That approach was subsequently used and
developed in the same framework by E. Mitidieri, S. Pokhozhaev and
many others, almost none of which cite the original research.

For a survey of the literature on the asymptotic behaviour and
blow-up of  solutions to the Cauchy problem for nonlinear parabolic
equations or inequalities we refer to [1], [3], [5], [14],  [16] and
[17].

\section{ Proofs}
{\it Proof of Theorem 1.} Let $n\geq 1$, $\alpha>0$ and  $1<q\leq 1
+\frac \alpha n$, let $\mathcal L$ be a differential operator
defined  by (2), the coefficients of which  satisfy the condition
(8) with the given $\alpha$ and some $c>0$,   and  let $(u,v)$ be an
entire weak solution of the inequality (1) in $\mathbb S$ such that
$u\geq v$. By the well-known inequality
$$(|u|^{q-1}u-|v|^{q-1}v)(u-v)\geq 2^{1-q} |u-v|^{q+1}$$ which holds
for any $q\geq 1$ and any  $u, v\in {\Bbb R}^1$, see, e.g., [6], we
obtain from (9) the relation  \begin{eqnarray} \int\limits_{\mathbb
S}\left [(u-v)_t\varphi +\sum_{i,j=1}^n a_{ij}(t,x)\frac{\partial
{\varphi}}{\partial x_i} \frac{\partial (u-v)}{\partial
x_j}\right]dtdx \geq 2^{1-q}\int\limits_{\mathbb S}(u-v)^q \varphi
dtdx
\end{eqnarray} which
holds for every  function $\varphi \in {C}{^\infty} (\mathbb S)$
with compact support in $\mathbb S$. Let $\tau>0$, $R>1$ and $T>0$
be real numbers. Let $\eta\! :[0,+\infty)\to [0,1]$ be a
$C^\infty$-function which has the non-negative  derivative $\eta'$
and equals 0 on the interval $[0, \tau]$ and 1 on the interval
$[2\tau, +\infty)$, and let $\zeta\! :[0,+\infty)\times {\Bbb R}^n
\to [0,1]$ be a $C^\infty$-function which equals 1 on $[0,
T/2]\times \overline {B (R/2)}$ and 0 on $\{[0, +\infty) \times
{\Bbb R}^n\} \setminus \{[0,T]\times \overline {B(R)}\}$, where
$B(R)$ is the ball in ${\Bbb R}^n$ centered at the origin of ${\Bbb
R}^n$ with radius $R$. Let
\begin{eqnarray}\varphi (t,x)= (w(t,x)+\varepsilon )^{-\nu}\zeta
^s(t,x)\eta^2(t),\nonumber\end{eqnarray} where
$w(t,x)=u(t,x)-v(t,x)$, $\varepsilon>0$ and the positive constants
$s>1$ and $1>\nu>0$ will be chosen below. Substituting the function
$\varphi$ in  (14) and then integrating by parts there   we obtain
\begin{eqnarray}
 -\frac s{1-\nu}
\int\limits_0^T\int\limits_{B(R)}(w+\varepsilon)^{1-\nu}
\zeta_{t}\zeta^{s-1}\eta^2dtdx -\frac 2{1-\nu}
\int\limits_0^T\int\limits_{B(R)}(w+\varepsilon)^{1-\nu} \zeta^s\eta' \eta
dtdx\nonumber\\ -\nu \int\limits_0^T\int\limits_{B(R)} \sum\limits_{i,j=1}^n a_{ij}(t,x)
\frac{\partial w}{\partial x_i} \frac{\partial w}{\partial x_j}
(w+\varepsilon )^{-\nu -1}\zeta^s \eta^2 dtdx\nonumber\\ +
s\int\limits_0^T\int\limits_{B(R)} \sum\limits_{i,j=1}^n a_{ij}(t,x)
\frac{\partial \zeta}{\partial x_i} \frac{\partial w}{\partial x_j} (w+\varepsilon)^{-\nu}\zeta^{s-1} \eta^2 dtdx
\nonumber \\ \equiv I_1+I_2+I_3 +I_4
\geq \int\limits_0^T\int\limits_{B(R)}w^q (w+\varepsilon )^{-\nu}\zeta^s \eta^2
dtdx.
\end{eqnarray}
In (15), first observe  that $I_3$   is  non-positive and then
estimate $I_4$ in terms of $I_3$. Namely, since
\begin{eqnarray}
|I_4|= \left| s\int\limits_0^T\int\limits_{B(R)} \sum\limits_{i,j=1}^n a_{ij}(t,x)
\frac{\partial \zeta}{\partial x_i} \frac{\partial w}{\partial x_j}
(w+\varepsilon )^{-\nu}\zeta^{s-1} \eta^2 dtdx \right|\nonumber \\
\leq \int\limits_0^T\int\limits_{B(R)} s \left(\sum\limits_{i,j=1}^n
a_{ij}(t,x) \frac{\partial w}{\partial x_i} \frac{\partial w}{\partial x_j}\right)^{\frac12} \left(\sum\limits_{i,j=1}^n
a_{ij}(t,x) \frac{\partial \zeta}{\partial x_i} \frac{\partial
\zeta}{\partial x_j}\right)^{\frac12}\nonumber \\
\cdot(w+\varepsilon )^{-\nu}\zeta^{s-1} \eta^2 dtdx, \end{eqnarray}
we estimate, first, the right-hand side of (16)  by using Young's
inequality
\begin{eqnarray}
AB\leq \rho A^2+ \rho^{-1} B^2, \nonumber\end{eqnarray} with
$\rho=\frac {\nu}2$,
$$A=\left(\sum\limits_{i,j=1}^n a_{ij}(t,x) \frac{\partial w}{\partial
x_i} \frac{\partial w}{\partial x_j}\right)^{\frac12}
(w+\varepsilon)^{-\frac{1+\nu}2}\zeta^{\frac s2}\eta$$ and $$
B=s\left(\sum\limits_{i,j=1}^n a_{ij}(t,x) \frac{\partial
\zeta}{\partial x_i} \frac{\partial \zeta}{\partial
x_j}\right)^{\frac12} (w+\varepsilon )^{\frac {1-\nu}2} \zeta^{\frac
s2 -1}\eta.$$ As a result, we arrive at
\begin{eqnarray}
|I_4|\leq \frac {\nu}2  \int\limits_0^T\int\limits_{B(R)} \sum\limits_{i,j=1}^n
a_{ij}(t,x) \frac{\partial w}{\partial x_i} \frac{\partial
w}{\partial x_j} (w+\varepsilon)^{-\nu-1}\zeta^s\eta^2 dtdx\nonumber\\
+\int\limits_0^T\int\limits_{B(R)}\frac {2s^2}\nu \sum\limits_{i,j=1}^n
a_{ij}(t,x) \frac{\partial \zeta}{\partial x_i} \frac{\partial
\zeta}{\partial x_j} (w+\varepsilon)^{1-\nu}\zeta^{s-2}\eta^2dtdx.
\end{eqnarray}
Further, since $I_2$ in (15) is also non-positive, the inequality
\begin{eqnarray}
\int\limits_0^T\int\limits_{B(R)}\frac s{1-\nu} (w+\varepsilon)^{1-\nu}
|\zeta_{t}|\zeta^{s-1}\eta^2dtdx\nonumber\\ +
\int\limits_0^T\int\limits_{B(R)}\frac {2s^2}\nu \sum\limits_{i,j=1}^n a_{ij}(t,x)
\frac{\partial \zeta}{\partial x_i} \frac{\partial \zeta}{\partial
x_j} (w+\varepsilon)^{1-\nu}\zeta^{s-2}\eta^2dtdx\nonumber \\ \geq
\int\limits_0^T\int\limits_{B(R)} w^q (w+\varepsilon)^{-\nu}\zeta^s \eta^2
dtdx\nonumber \\
+\frac{\nu}2\int\limits_0^T\int\limits_{B(R)}\sum\limits_{i,j=1}^n
a_{ij}(t,x) \frac{\partial w}{\partial x_i}\frac{\partial
w}{\partial x_j} (w+\varepsilon)^{-\nu-1}\zeta^s\eta^2 dtdx
\end{eqnarray} easily follows from (15) and (17). Estimating both
integrands on the left-hand side of (18) by Young's inequality
\begin{eqnarray}
AB\leq \rho A^{\frac \beta{\beta-1}}+\rho^{1-\beta} B^\beta,\nonumber
\end{eqnarray} respectively,  with
$\rho=\frac 14$, $\beta=\frac{q-\nu}{q-1}$,
$$A=(w+\varepsilon)^{1-\nu}{\zeta}^{\frac
{s(1-\nu)}{q-\nu}}\eta^{\frac{2(1-\nu)}{q-\nu}},$$
$$B=\frac s{1-\nu}
|\zeta_t| {\zeta}^{{\frac
{s(q-1)}{q-\nu}}-1}\eta^{\frac{2(q-1)}{q-\nu}}$$  and with $\rho=\frac
{1}4$, $\beta=\frac{q-\nu}{q-1}$,
$$A=(w+\varepsilon)^{1-\nu}{\zeta}^{\frac
{s(1-\nu)}{q-\nu}}\eta^{\frac{2(1-\nu)}{q-\nu}},$$$$ B=\frac
{2s^2}\nu {\zeta}^{{\frac
{s(q-1)}{q-\nu}}-2}\eta^{\frac{2(q-1)}{q-\nu}} \sum\limits_{i,j=1}^n
a_{ij}(t,x) \frac{\partial \zeta}{\partial x_i} \frac{\partial
\zeta}{\partial x_j},
$$  we obtain
\begin{eqnarray}
\frac {1}4 \int\limits_0^T\int\limits_{B(R)}(w+\varepsilon)^{q-\nu}\zeta^s \eta^2
dtdx +c_1\int\limits_0^T\int\limits_{B(R)}|\zeta_t|^{\frac {q-\nu}{q-1}}
\zeta^{s-{\frac{q-\nu}{q-1}}}\eta^2dtdx\nonumber\\ +\frac {1}4\int\limits_0^T\int\limits_{B(R)}(w+\varepsilon)^{q-\nu}\zeta^s \eta^2
dtdx\nonumber \\ + c_2 \int\limits_0^T\int\limits_{B(R)}
\left(\sum\limits_{i,j=1}^n a_{ij}(t,x) \frac{\partial
\zeta}{\partial x_i} \frac{\partial \zeta}{\partial x_j}\right)^
{\frac {q-\nu}{q-1}}
\zeta^{s-{\frac{2(q-\nu)}{q-1}}}\eta^2dtdx\nonumber \\ \geq
\int\limits_0^T\int\limits_{B(R)}w^q (w+\varepsilon)^{-\nu}\zeta^s \eta^2 dtdx
\nonumber \\+ \frac\nu{2}\int\limits_0^T\int\limits_{B(R)}\sum\limits_{i,j=1}^n
a_{ij}(t,x) \frac{\partial w}{\partial x_i}\frac{\partial
w}{\partial x_j} (w+\varepsilon)^{-\nu-1}\zeta^s \eta^2 dtdx.
\end{eqnarray}
Here and what follows, we use the symbols  $c_i$, $i=1, 2, \ldots$,  to denote constants depending possibly  on $c$, $n$, $q$, $s$, $\alpha$
and  $\nu$,   but not on $\varepsilon$, $\tau$  and $R$.

At this point,  using the inequality (19),  we obtain an upper bound
on the integral
\begin{eqnarray}
\int\limits_0^T\int\limits_{B(R)}w^q\zeta^s \eta^2 dtdx. \nonumber\end{eqnarray} To
this end, we substitute
$$\varphi (t,x) = \zeta^s(t,x) \eta^2(t)$$
in the inequality (9) and  then after the integration by parts there
we have
\begin{eqnarray}
-s\int\limits_0^T\int\limits_{B(R)}w\zeta_t\zeta^{s-1}\eta^2dtdx -2
\int\limits_0^T\int\limits_{B(R)}w\zeta^s\eta' \eta dtdx \nonumber
\\ +s\int\limits_0^T\int\limits_{B(R)}\sum\limits_{i,j=1}^n a_{ij}(t,x)
\frac{\partial \zeta}{\partial x_i}\frac{\partial w}{\partial x_j}
\zeta^{s-1}\eta^2 dtdx \geq
\int\limits_0^T\int\limits_{B(R)}w^q\zeta^s \eta^2 dtdx.
\end{eqnarray} As before,  it is easy to see that the second term on
the left-hand side of (20) is non-positive and thus (20) yields
\begin{eqnarray}
s\int\limits_0^T\int\limits_{B(R)}w|\zeta_t|\zeta^{s-1}\eta^2dtdx
+s\int\limits_0^T\int\limits_{B(R)}\sum\limits_{i,j=1}^n a_{ij}(t,x)
\frac{\partial \zeta}{\partial x_i}\frac{\partial w}{\partial x_j}
\zeta^{s-1}\eta^2 dtdx\nonumber \\ \geq
\int\limits_0^T\int\limits_{B(R)}w^q\zeta^s \eta^2
dtdx.\end{eqnarray} Estimating now the first integral on the
left-hand side of (21) by H\"{o}lder's inequality, we arrive at
\begin{eqnarray}
s\left(\int\limits_{T/2}^T\int\limits_{B(R)}w^q\zeta^s \eta^2
dtdx \right)^{\frac{1}{q}}
\left(\int\limits_0^T\int\limits_{B(R)}|\zeta_t|^{\frac q {q-1}}
\zeta^{s-\frac q{q-1}} \eta^2 dtdx\right)^{\frac{q-1}{q}} \nonumber
\\+s\int\limits_0^T\int\limits_{B(R)}\sum\limits_{i,j=1}^n a_{ij}(t,x)
\frac{\partial \zeta}{\partial x_i}\frac{\partial w}{\partial x_j}
\zeta^{s-1}\eta^2 dtdx \geq \int\limits_0^T\int\limits_{B(R)}w^q\zeta^s \eta^2
dtdx. \end{eqnarray} Further, since
\begin{eqnarray}
\left|\sum\limits_{i,j=1}^n a_{ij}(t,x) \frac{\partial \zeta}{\partial
x_i} \frac{\partial w}{\partial x_j}\right|
\nonumber \\
\leq \left(\sum\limits_{i,j=1}^n a_{ij}(t,x) \frac{\partial
w}{\partial x_i} \frac{\partial w}{\partial x_j}\right)^{\frac12}
\left(\sum\limits_{i,j=1}^n a_{ij}(t,x) \frac{\partial
\zeta}{\partial x_i} \frac{\partial \zeta}{\partial
x_j}\right)^{\frac12},\nonumber
\end{eqnarray} we have
\begin{eqnarray}
\int\limits_0^T\int\limits_{B(R)}\sum\limits_{i,j=1}^n a_{ij}(t,x) \frac{\partial
\zeta}{\partial x_i}\frac{\partial w}{\partial x_j}
\zeta^{s-1}\eta^2 dtdx\nonumber  \\ \leq \int\limits_0^T\int\limits_{B(R)}
\left(\sum\limits_{i,j=1}^n a_{ij}(t,x) \frac{\partial \zeta}{\partial
x_i} \frac{\partial \zeta}{\partial x_j}\right)^{\frac 12}
\left(\sum\limits_{i,j=1}^n a_{ij}(t,x) \frac{\partial
w}{\partial x_i} \frac{\partial w}{\partial x_j}\right)^
{\frac 12} \zeta^{s-1}\eta^2 dtdx.
\end{eqnarray}
Estimating the right-hand side of (23) by H\"{o}lder's inequality,
we obtain the relation
\begin{eqnarray}
\int\limits_0^T\int\limits_{B(R)} \sum\limits_{i,j=1}^n a_{ij}(t,x) \frac{\partial
\zeta}{\partial x_i} \frac{\partial w}{\partial x_j}
\zeta^{s-1}\eta^2 dtdx\nonumber \\ \leq \left(\int\limits_0^T\int\limits_{B(R)}
\sum\limits_{i,j=1}^n a_{ij}(t,x) \frac{\partial \zeta}{\partial
x_i} \frac{\partial \zeta}{\partial x_j} (w+\varepsilon )^{1+\nu}
\zeta^{s-2}\eta^2 dtdx\right)^{1/2}\nonumber \\
\times \left(\int\limits_0^T\int\limits_{B(R)}\sum\limits_{i,j=1}^n a_{ij}(t,x)
\frac{\partial w}{\partial x_i}\frac{\partial w}{\partial x_j}
(w+\varepsilon) ^{-\nu-1}\zeta^s \eta^2 dtdx \right)^{1/2}
\end{eqnarray} which holds for any  $\varepsilon >0$ and any  $\nu\in (0,1)$.
Further, it is  easy to see that the inequality
\begin{eqnarray}
\int\limits_0^T\int\limits_{B(R)}\sum\limits_{i,j=1}^n a_{ij}(t,x) \frac{\partial
\zeta}{\partial x_i} \frac{\partial \zeta}{\partial x_j}
(w+\varepsilon)^{1+\nu}  \zeta^{s-2} \eta^2 dtdx\nonumber
\\ \leq \left (\int\limits_0^T\int\limits_{B(R)}
\left(\sum\limits_{i,j=1}^n a_{ij}(t,x) \frac{\partial
\zeta}{\partial x_i} \frac{\partial \zeta}{\partial x_j}\right)^
{\frac {d}{d-1}} \zeta^{s-\frac{2d}{d-1}}\eta^2 dtdx\right)^{\frac
{d-1} d}\nonumber \\\times  \left(\int\limits_0^T\int\limits_{B(R)\setminus B(R/2)}
(w+\varepsilon)^{d(1+\nu)}\zeta^s \eta^2 dtdx \right)^{\frac 1 d}
\end{eqnarray}
holds for any $d>1$. In (25), choosing  for any  sufficiently small
$\nu\in(0, 1)$ the parameter $d$ such that $d(1+\nu)=q$, we obtain
from (24) and (25)  the relation
\begin{eqnarray}
\int\limits_0^T\int\limits_{B(R)}\sum\limits_{i,j=1}^n a_{ij}(t,x) \frac{\partial
\zeta}{\partial x_i}\frac{\partial w}{\partial x_j}
\zeta^{s-1}\eta^2 dtdx
\nonumber
\\ \leq
\left (\int\limits_0^T\int\limits_{B(R)} \left(\sum\limits_{i,j=1}^n
a_{ij}(t,x) \frac{\partial \zeta}{\partial x_i} \frac{\partial
\zeta}{\partial x_j}\right)^ {\frac {d}{d-1}} \zeta^
{s-\frac{2d}{d-1}}\eta^2 dtdx\right)^{\frac {d-1} {2d}}\nonumber \\
\times  \left(\int\limits_0^T\int\limits_{B(R)\setminus B(R/2)}
(w+\varepsilon)^{q} \zeta^s
\eta^2 dtdx \right)^{\frac 1 {2d}}\nonumber \\
\times\left(\int\limits_0^T\int\limits_{B(R)}\sum\limits_{i,j=1}^n a_{ij}(t,x)
\frac{\partial w}{\partial x_i}\frac{\partial w}{\partial x_j}
(w+\varepsilon)^{-\nu-1}\zeta^s \eta^2 dtdx \right)^{1/2}
\end{eqnarray}
which holds for any $\varepsilon >0$ and any sufficiently small
$\nu\in (0, 1)$.  In (26), estimating the last term on the
right-hand side  by virtue of (19), we have
\begin{eqnarray}
\int\limits_0^T\int\limits_{B(R)}\sum\limits_{i,j=1}^n a_{ij}(t,x) \frac{\partial
\zeta}{\partial x_i}\frac{\partial w}{\partial x_j}
\zeta^{s-1}\eta^2 dtdx \nonumber \\ \leq \left(\int\limits_0^T\int\limits_{B(R)}
\left(\sum\limits_{i,j=1}^n a_{ij}(t,x) \frac{\partial
\zeta}{\partial x_i} \frac{\partial \zeta}{\partial
x_j}\right)^{\frac {d}{d-1}} \zeta^{s-\frac{2d}{d-1}}
\eta^2 dtdx\right)^{\frac {d-1}{2d}} \nonumber \\ \times
\left(\int\limits_0^T\int\limits_{B(R)\setminus B(R/2)} (w+\varepsilon)^{q} \zeta^s
\eta^2 dtdx \right)^{\frac 1 {2d}} \left( \frac {1}{\nu}
\int\limits_0^T\int\limits_{B(R)}(w+\varepsilon)^{q-\nu}\zeta^s \eta^2 dtdx\right.
\nonumber \\ -\frac {2}{\nu} \int\limits_0^T\int\limits_{B(R)}w^q(w+\varepsilon)^
{-\nu}\zeta^s\eta^2 dtdx + c_3\int\limits_0^T\int\limits_{B(R)}|\zeta_t|^{\frac
{q-\nu}{q-1}}
\zeta^{s-{\frac{q-\nu}{q-1}}}\eta^2 dtdx\nonumber \\
+\left.c_4\int\limits_0^T\int\limits_{B(R)} \left(\sum\limits_{i,j=1}^n a_{ij}(t,x)
\frac{\partial \zeta}{\partial x_i} \frac{\partial \zeta}{\partial
x_j}\right)^{\frac {q-\nu}{q-1}} \zeta^{s-{\frac
{2(q-\nu)}{q-1}}}\eta^2 dtdx\right)^{\frac{1}{2}}.
\end {eqnarray}
In (27), passing to the limit as $\varepsilon\to 0$ as justified by
Lebesgue's theorem (see, e.g., [9, p. 303]), we obtain for any
sufficiently large $s$ the inequality
\begin{eqnarray}
\int\limits_0^T\int\limits_{B(R)}\sum\limits_{i,j=1}^n a_{ij}(t,x) \frac{\partial
\zeta}{\partial x_i}\frac{\partial w}{\partial x_j}
\zeta^{s-1}\eta^2 dtdx\nonumber \\ \leq  \left
(\int\limits_0^T\int\limits_{B(R)} \left(\sum\limits_{i,j=1}^n a_{ij}(t,x)
\frac{\partial \zeta}{\partial x_i} \frac{\partial \zeta}{\partial
x_j}\right)^
{\frac {d}{d-1}}\eta^2 dtdx\right)^{\frac {d-1}{2d}} \nonumber \\
\times \left(\int\limits_0^T\int\limits_{B(R)\setminus B(R/2)}w^q \zeta^s \eta^2 dt dx \right)^{\frac 1 {2d}}\nonumber \\
\times \left(c_3\int\limits_0^T\int\limits_{B(R)}|\zeta_t|^{\frac {q-\nu}{q-1}}\eta^2 dt dx
+c_4\int\limits_0^T\int\limits_{B(R)} \left(\sum\limits_{i,j=1}^n a_{ij}(t,x)
\frac{\partial \zeta}{\partial x_i} \frac{\partial \zeta}{\partial
x_j}\right)^ {\frac {q-\nu}{q-1}}\eta^2 dtdx\right)^{\frac 12}.
\end{eqnarray}
In turn, (22) and (28) yield the inequality
\begin{eqnarray}
\int\limits_0^{T}\int\limits_{B(R)}w^q\zeta^s \eta^2 dtdx\nonumber \\ \leq
c_5\left(\int\limits_{T/2}^T\int\limits_{B(R)}w^q\zeta^s\eta^2
dtdx\right)^{\frac{1}{q}}
\left(\int\limits_0^T\int\limits_{B(R)}|\zeta_t|^{\frac{q}{q-1}}\eta^2dtdx\right)^{\frac{q-1}{q}}
\nonumber \\
+c_5\left (\int\limits_0^T\int\limits_{B(R)} \left(\sum\limits_{i,j=1}^n
a_{ij}(t,x) \frac{\partial \zeta}{\partial x_i} \frac{\partial
\zeta}{\partial x_j}\right)^ {\frac {d}{d-1}} \eta^2 dtdx\right)^{\frac
{d-1} {2d}}\nonumber \\ \times \left(\int\limits_0^T\int\limits_{B(R)\setminus B(R/2)} w^q \zeta^s \eta^2 dtdx
\right)^{\frac 1 {2d}}\nonumber \\
\times\left(\int\limits_0^T\int\limits_{B(R)}|\zeta_t|^{\frac{q-\nu}{q-1}}  \eta^2dtdx
+\int\limits_0^T\int\limits_{B(R)} \left(\sum\limits_{i,j=1}^n a_{ij}(t,x)
\frac{\partial \zeta}{\partial x_i} \frac{\partial \zeta}{\partial
x_j}\right) ^{\frac{q-\nu}{q-1}} \eta^2  dt dx  \right)^{\frac 12}
\end{eqnarray} which holds for any sufficiently large $s$  and any sufficiently small $\nu\in (0,1)$.
Further, by the condition (3) on the coefficients of the operator $\mathcal L$, we obtain
from (29) the inequality
\begin{eqnarray}
\int\limits_0^T\int\limits_{B(R)}w^q\zeta^s \eta^2 dtdx\nonumber \\ \leq
c_5\left(\int\limits_{T/2}^T\int\limits_{B(R)}w^q\zeta^s\eta^2
dtdx\right)^{\frac{1}{q}}
\left(\int\limits_0^T\int\limits_{B(R)}|\zeta_t|^{\frac{q}{q-1}}\eta^2 dtdx\right)^{\frac{q-1}{q}}
\nonumber \\
+c_5\left({\mathcal A}(T,R)\right)^\frac 12 \left(\int\limits_0^T\int\limits_{B(R)} |\nabla \zeta|^{\frac {2d}{d-1}}\eta^2 dtdx
\right)^{\frac{d-1}{2d}}
\left(\int\limits_0^T\int\limits_{B(R)\setminus
B(R/2)} w^q \zeta^s \eta^2 dtdx
\right)^{\frac 1 {2d}}\nonumber \\
\times\left(\int\limits_0^T\int\limits_{B(R)}|\zeta_t|^{\frac{q-\nu}{q-1}} \eta^2 dtdx
+
\left({\mathcal A}(T, R)\right)^\frac{q-\nu}{q-1}
\int\limits_0^T\int\limits_{B(R)}( |\nabla \zeta|^{\frac{2(q-\nu)}{q-1}} \eta^2 dtdx \right)^{\frac 12},
\nonumber\end{eqnarray} where
\begin{eqnarray}{\mathcal A}(T, R)={\mathrm {ess}\sup}_{(t,x)\in
(0,T )\times \{B(R)\setminus B(R/2)\}}A(t,x), \nonumber \end{eqnarray}
which in turn by the condition (8) yields
\begin{eqnarray}
\int\limits_0^T\int\limits_{B(R)}w^q\zeta^s \eta^2 dtdx\nonumber \\ \leq
c_6\left(\int\limits_{T/2}^T\int\limits_{B(R)}w^q\zeta^s\eta^2
dtdx\right)^{\frac{1}{q}}
\left(\int\limits_0^T\int\limits_{B(R)}|\zeta_t|^{\frac{q}{q-1}}\eta^2  dtdx\right)^{\frac{q-1}{q}}
\nonumber \\
+c_6R^{\frac {2-\alpha}{2}}\left(\int\limits_0^T\int\limits_{B(R)} |\nabla \zeta|^{\frac {2d}{d-1}}\eta^2 dtdx
\right)^{\frac{d-1}{2d}}
\left(\int\limits_0^T\int\limits_{B(R)\setminus
B(R/2)} w^q \zeta^s \eta^2 dtdx
\right)^{\frac 1 {2d}}\nonumber \\
\times\left(\int\limits_0^T\int\limits_{B(R)}|\zeta_t|^{\frac{q-\nu}{q-1}} \eta^2  dtdx
+
R^{\frac{(2-\alpha)(q-\nu)}{(q-1)}}
\int\limits_0^T\int\limits_{B(R)} |\nabla \zeta|^{\frac{2(q-\nu)}{q-1}} \eta^2 dtdx \right)^{\frac 12}.
\end{eqnarray}
Now, for arbitrary $(t,x)\in \mathbb S$, $R>1$ and $T>0$, we choose in (30)
the function $\zeta(t, x)$ in the form
\begin{equation}
\zeta
(t, x)= \psi\left(\frac tT\right)  \psi\left(\frac {2|x|^2}{R^2}\right),   \end{equation}
where
$\psi: [0, +\infty)\to [0,1]$ is a $C^\infty$-function which equals 1 on
$[0, 1/2]$ and 0 on $[1, +\infty)$ and such that the inequalities
\begin{eqnarray}
|\zeta_t|\leq {c_7}T^{-1}\qquad \hbox{and}\qquad  |\nabla_x\zeta|\leq
{c_7}R^{-1}\end {eqnarray}
hold.
Note that
it is always possible to find such a function $\zeta$. Indeed, this
can be easily verified by direct calculation of the corresponding
derivatives  of the function $\zeta$  defined
by (31). Also, in what follows we let
\begin{eqnarray}T=R^\alpha.\end{eqnarray}
Since $|\eta|\leq 1$, by (32) and (33),  we have  from  (30) the inequality
\begin{eqnarray}
\int\limits_0^{T}\int\limits_{B(R)}w^q\zeta^s \eta^2 dtdx \leq
c_8\left(R^{n+\alpha-\frac {\alpha q}
{q-1}}\right)^{\frac{q-1}q}\left(\int\limits_{T/2}^T\int\limits_{B(R)}w^q\zeta^s\eta^2dtdx\right)^{\frac{1}{q}}\nonumber \\  +
c_8
R^{\frac{2-\alpha}2}
\left(R^{n+\alpha-\frac
{2d}{d-1}}\right)^{\frac{d-1}{2d}} \left(R^{n+\alpha-\frac {\alpha (q-\nu)} {q-1}} +
R^{\frac{(2-\alpha)(q-\nu)}{q-1}}
R^{n+\alpha-\frac {2(q-\nu)} {q-1}}
\right)^{\frac 12}\nonumber
\\ \times
\left(\int\limits_0^{T}\int\limits_{B(R)\setminus
B(R/2)} w^q \zeta^s \eta^2 dtdx \right)^{\frac 1 {2d}}.
\end{eqnarray}
Making simple calculations in (34) we obtain
\begin{eqnarray}
\int\limits_0^T\int\limits_{B(R)}w^q\zeta^s \eta^2 dtdx \leq
c_8\left(R^{n+\alpha-\frac {\alpha q}
{q-1}}\right)^{\frac{q-1}q}\left(\int\limits_{T/2}^T\int\limits_{B(R)}w^q\zeta^s\eta^2dtdx\right)^{\frac{1}{q}}\nonumber \\  +
c_8
\left(R^{n+\alpha-\frac
{d\alpha}{d-1}}\right)^{\frac{d-1}{2d}} \left(R^{n+\alpha-\frac {\alpha (q-\nu)} {q-1}}
\right)^{\frac 12}
\left(\int\limits_0^{T}\int\limits_{B(R)\setminus
B(R/2)} w^q \zeta^s \eta^2 dtdx \right)^{\frac 1 {2d}}.
\end{eqnarray}
In turn, since $d(1+\nu)=q$, i.e.,
$\frac{2d}{d-1}=\frac{2q}{q-1-\nu}$, the relation (35) implies the inequality
\begin{eqnarray}
\int\limits_0^T\int\limits_{B(R)}w^q\zeta^s \eta^2 dtdx\nonumber \\ \leq
c_8R^{\frac{n}{q-1}[q-1-\frac{\alpha}{n}]}\left(\int\limits_{T/2}^T\int\limits_{B(R)}w^q\zeta^s\eta^2 dtdx\right)^{\frac{1}{q}}\nonumber \\
+c_8R^{\frac {n(2q-1-\nu)}{2 q(q-1)}\left[q-1-\frac {\alpha}n\right]}
\left(\int\limits_0^T\int\limits_{B(R)\setminus B(R/2)} w^q \zeta^s \eta^2 dtdx
\right)^{\frac 1 {2d}} \end{eqnarray}
which holds
for any sufficiently small $\nu\in (0,1)$.
Now, since $q>1$, $d>1$, and for any $\nu\in (0,1)$ both quantities
\begin{eqnarray}
\frac {n}{q-1} \qquad {\rm and} \qquad \frac
{n(2q-1-\nu)}{2q(q-1)} \nonumber\end{eqnarray} are  positive,
it follows from (36), where we remind that  $T=R^\alpha$,  by passing $R\to+\infty$   that the relation
\begin{eqnarray}\int_{\mathbb S} w^q \eta^2dtdx=0\end{eqnarray} holds
for \begin{eqnarray}1<q<1+\frac {\alpha}n.\nonumber\end{eqnarray} We
show now that (37) also holds for \begin{eqnarray}q =1+ \frac
{\alpha}n.\end{eqnarray} Indeed, since $q>1$ and $d>1$, by (38) we
have from (36) the estimate
\begin{eqnarray}\int_{\mathbb S} w^q \eta^2dtdx<+\infty.
\end{eqnarray}
In turn,  by Fubini's theorem (see, e.g., [9, p. 317]), we obtain
from (39) the relations
\begin{eqnarray}
\int\limits_{T_k/2}^{T_k}\int\limits_{{\Bbb R}^n} w^q \eta^2 dtdx\to 0
\end{eqnarray} and
\begin{eqnarray}
\int\limits_0^{+\infty}\int\limits_{B(R_k)\setminus B(R_k/2)}  w^q \eta^2 dtdx\to 0
\end{eqnarray}
which hold for any sequences $R_k$ and $T_k$ such that  $R_k\to
+\infty$ and $T_k\to +\infty$. On the other hand, from (36) we have
the  inequality
\begin{eqnarray}
\int\limits_0^{T/2}\int\limits_{B(R/2)}w^q \eta^2dtdx \leq c_{8}R^{\frac
n{q-1}\left[q-1-\frac
{\alpha}n\right]}\left(\int\limits_{T/2}^T\int\limits_{{\mathbb R}^n}w^q\eta^2dtdx\right)^{\frac{1}{q}}\nonumber \\
+c_{8}R^{\frac {n(2q-1-\nu)}{2q(q-1)}\left[q-1-\frac
{\alpha}n\right]} \left(\int\limits_0^{+\infty}\int\limits_{B(R)\setminus
B(R/2)} w^q \eta^2 dtdx \right)^{\frac 1 {2d}}
\end{eqnarray} which,  together with (40) and (41), where we choose $T=R^\alpha=T_k=(R_k)^\alpha$,  implies  the relation
\begin{eqnarray}
\int\limits_0^{(R_k/2)^\alpha}\int\limits_{B(R_k/2)}w^q \eta^2dtdx \to 0
\end{eqnarray} which holds for any  sequence $R_k\to +\infty$. In
turn, (43) yields the relation (37) for $q$ given by (38). Thus, we
prove that the relation (37) holds for any \begin{eqnarray}1<q \leq
1+\frac {\alpha}n,\nonumber\end{eqnarray} where $\eta\! :
[0,+\infty)\to [0,1]$ is a $C^\infty$-function which equals 1 on the
interval $[2\tau, +\infty)$. In (37), passing to the limit as
$\tau\to 0$, we obtain that $u(t,x)=v(t,x)$  in $\mathbb S$.

\vspace{10mm}

\newpage

\noindent \textbf{Authors' addresses:}

\vspace{5 mm}

\noindent Vasilii V. Kurta

\noindent Mathematical Reviews

\noindent 416 Fourth Street, P.O. Box 8604

\noindent Ann Arbor, Michigan 48107-8604, USA

\noindent \textbf {e-mail:} vkurta@umich.edu, vvk@ams.org
\end{document}